\documentclass[11pt, a4paper]{amsart}
\usepackage[T1]{fontenc}
\usepackage[utf8]{inputenc}
\usepackage{lmodern,amsmath,amssymb,microtype}
\usepackage[a4paper,left=28mm,right=28mm,top=23mm,bottom=23mm]{geometry}
\usepackage{xcolor,etoolbox}
\definecolor{links}{RGB}{26,66,110}
\definecolor{citations}{RGB}{18,84,75}
\usepackage[colorlinks=true,linkcolor=links,citecolor=citations,urlcolor=links,
 pdfdisplaydoctitle=true,pdfauthor={Leokadia Białas-Cież, Rafał Pierzchała},
 pdftitle={Holder continuity of Green functions for Markov sets}]{hyperref}
\newtheorem{theorem}{Theorem}[section]
\newtheorem{lemma}[theorem]{Lemma}
\newtheorem{corollary}[theorem]{Corollary}
\theoremstyle{definition}
\newtheorem{definition}[theorem]{Definition}
\theoremstyle{plain}
\newcommand{\C}{\mathbb C}
\newcommand{\Pn}{\mathcal P}
\makeatletter

\newcommand{\normindex}[2]{\lower.17ex\hbox{$\m@th#1#2$}}
\newcommand{\norm}[2]{{\lVert#1\rVert}_{\mathpalette\normindex{#2}}}
\makeatother
\newcommand{\dist}{\operatorname{dist}}

\makeatletter

\patchcmd{\@maketitle}{\global\topskip42\p@}{\global\topskip34\p@}{}{}
\patchcmd{\@starttoc}{32\p@\@plus14\p@}{12\p@\@plus2\p@}{}{}
\patchcmd{\@starttoc}{\linespacing\@plus\linespacing}{.8\linespacing\@plus2\p@}{}{}
\renewcommand{\tocsection}[3]{\makebox[1.9em][l]{\@ifnotempty{#2}{#2.}}#3}

\renewcommand{\section}{\@startsection{section}{1}{\z@}%
  {1.55\linespacing \@plus .35\linespacing \@minus .16\linespacing}%
  {.80\linespacing \@plus .12\linespacing}{\normalfont\scshape\centering}}
\renewcommand{\subsection}{\@startsection{subsection}{2}{\z@}%
  {1.10\linespacing \@plus .22\linespacing \@minus .12\linespacing}%
  {.48\linespacing}{\normalfont\bfseries}}
\renewcommand{\l@section}{\@tocline{1}{0pt}{1pc}{0pt}{}}
\makeatother
\title[Markov's inequality and the Green's function]{H\"older continuity of Green functions \texorpdfstring{\\}{ }for Markov sets}

\author[L. Bialas-Ciez]{Leokadia Bialas-Ciez}
\address{Faculty of Mathematics and Computer Science, Jagiellonian University,
Łojasiewicza 6, 30-348 Kraków, Poland}
\email{\href{mailto:leokadia.bialas-ciez@uj.edu.pl}{leokadia.bialas-ciez@uj.edu.pl}}
\author[R. Pierzchała]{Rafał Pierzchała}
\address{Faculty of Mathematics and Computer Science, Jagiellonian University,
Łojasiewicza 6, 30-348 Kraków, Poland}
\email{\href{mailto:rafal.pierzchala@im.uj.edu.pl}{rafal.pierzchala@uj.edu.pl}}

\date{}
\begin{document}
\begin{abstract}
Let $K$ be a nonempty compact subset of $\mathbb C^N$. We prove that
its Green function (and consequently, the Siciak extremal function) is Hölder continuous if and only if the
first derivatives of polynomials satisfy Markov's inequality on $K$. The implication from Markov’s inequality to Hölder continuity of the associated Green function settles a longstanding open problem at the intersection of approximation theory and pluripotential theory. 
A Markov bound $An^m$ in degree $n$ gives every Hölder exponent below
$1/m$. In one complex variable it gives the exponent $1/m$ itself.
The proof compares the extremal functions of two neighborhoods of $K$.
For $N>1$ we use Hörmander's theorem on holomorphic extension from a
point. For $N=1$ the comparison follows from the maximum principle.
\end{abstract}
\maketitle

\section{The Markov and Hölder inequalities}\label{sec:introduction}

We ask how Markov's inequality on a compact set $K\subset\C^N$ controls
polynomial growth near $K$. Siciak's extremal function $V_K$ measures this growth.
If $K$ satisfies Markov's inequality with exponent~$1$, a direct calculation
shows that $V_K$ is Lipschitz continuous on $\C^N$
\cite{BBM,BC11}. See also Corollary~\ref{cor:linear}. The question is what
happens when the Markov exponent is greater than~$1$.

In 1990 Pleśniak posed two problems \cite[p.~116, Open Problems~$1^\circ$ and~$2^\circ$]{PlesExt}
about compact sets $K\subset\mathbb R^N$ satisfying Markov's inequality.
In our notation, he asked whether $V_K$ is continuous on $\C^N$
(equivalently, whether $K$ is $L$-regular) and, in particular, whether $K$ is nonpluripolar. He also asked for an example of such a set without the
Hölder continuity property (HCP).
In his Oberwolfach report \cite[p.~13]{Ples91} he returned to the
questions of $L$-regularity and positive logarithmic capacity in one variable.
Siciak asked about $L$-regularity and positive capacity for compact Markov
sets in $\mathbb R^N$ and $\C^N$ \cite[pp.~447--448]{Siciak95},
and about HCP in the plane \cite[p.~219]{Sic97plane}. Levenberg asked
whether compact Markov sets in $\C^N$ have HCP, are $L$-regular and,
for $N>1$, are nonpluripolar \cite[p.~135, Problem~6]{Levenberg}.
See also \cite{BB14,BB24,BBM,BC11,P22}.

Pleśniak \cite{PlesCantor} constructed an $L$-regular Cantor set on the real line for which Markov’s inequality fails. See also \cite[Example~5.1]{Kosek21}.
HCP and Markov's inequality for the middle-thirds
Cantor set was proved in \cite{BV93}, which answered another Pleśniak's problem. For certain classes of Cantor sets on the real line, the
equivalence with HCP was proved in \cite{BC95,Tot}. Every planar Markov set is nonpolar, as shown in \cite{BC98}.
$L$-regularity of compact Markov subsets of the real line and of compact
sets satisfying a local Markov inequality in the complex plane was proved
by Bialas-Ciez and Eggink \cite{BE08}.

We answer these questions by proving that every nonempty compact Markov
set in $\C^N$ has HCP. Hence it is $L$-regular and nonpluripolar.
In one variable it has positive logarithmic capacity. The example sought
in Pleśniak's second problem therefore does not exist.

Baran and Bialas-Ciez characterized HCP by the V.~Markov inequality
\cite[Definition~2.5 and Theorem~2.9]{BB14}. It bounds each partial
derivative of order~$k$ of a polynomial of degree~$n$ by
$M^kn^{mk}/(k!)^{m-1}$ times the norm of the polynomial, with fixed
$M>0$ and $m\ge1$. For $m>1$ and $k\ge2$, the factorial improves
the iterated first-derivative bound $M^kn^{mk}$. Our main theorem
assumes only the first-derivative inequality.

For the nonpluripolarity question in several variables, see also
\cite[after Theorem~1.12]{P21}. If a nonempty compact Markov set
$K\subset\mathbb R^N$ is definable in an o-minimal expansion of the real
field, then $K=\overline{\operatorname{int}_{\mathbb R^N}K}$
\cite[Theorem~1.15]{P22}. In particular, $K$ is nonpluripolar.

\subsection*{Definitions and notation}
Throughout the paper, $K$ is a nonempty compact subset of $\C^N$, where $N\ge1$.
\begin{itemize}
\item We use the Euclidean norm $|\cdot|$ on $\C^N$.
      The open ball with center $a$ and radius $s>0$ is denoted by $B(a,s)$.
\item For $r\ge0$, the closed $r$-neighborhood of $K$ is
      \[
        K_r=\bigl\{z\in\C^N : \dist(z,K)\le r\bigr\}.
      \]
\item For a continuous scalar- or vector-valued function $f$ on $K$, we write
      \[
        \norm f K=\max_{x\in K}|f(x)|.
      \]
\item We write $\mathbb N=\{1,2,\ldots\}$ and $\mathbb N_0=\mathbb N\cup\{0\}$.
      For $n\in\mathbb N_0$, $\Pn_n$ denotes the space of complex polynomials
      of total degree at most $n$. The number of variables is understood from the context.
\item Write $\partial_j=\partial/\partial z_j$. For $v\in\C^N$, the
      complex gradient and the directional derivative of $P$ are
      \[
        \nabla P=(\partial_1P,\ldots,\partial_NP),\qquad
        D_vP=\sum_{j=1}^N v_j\partial_jP.
      \]
\item For $\alpha\in\mathbb N_0^N$ and $h\in\C^N$, we use the notation
      \[
        D^\alpha=\partial_1^{\alpha_1}\cdots\partial_N^{\alpha_N},\qquad
        h^\alpha=h_1^{\alpha_1}\cdots h_N^{\alpha_N},\qquad
        |\alpha|=\sum_{j=1}^N\alpha_j,\qquad
        \alpha!=\prod_{j=1}^N\alpha_j!.
      \]
\end{itemize}

\begin{definition}[Markov's inequality]\label{def:markov}
We say that $K$ satisfies Markov's inequality with exponent $m\ge1$,
or that $K$ is a Markov set with exponent $m$, if there is a constant $A>0$ such that
\begin{equation}\label{eq:markov}
 \norm{\nabla P}{K}\le An^m\norm{P}{K}
 \qquad\bigl(n\in\mathbb N,\ P\in\Pn_n\bigr).
\end{equation}
\end{definition}

Clearly, \eqref{eq:markov} is equivalent to the existence of a constant $M>0$
such that, for the same $m$,
\[
 \norm{D^\alpha P}{K}\le M^{|\alpha|}n^{m|\alpha|}\norm P K
 \qquad\bigl(n\in\mathbb N,\ P\in\Pn_n,\ \alpha\in\mathbb N_0^N\bigr).
\]
The prototype is A.~A.~Markov's classical inequality on the interval $I=[-1,1]$.
For each polynomial $p$ in one variable of degree at most $n$,
\[
                       \norm{p'}I\le n^2\norm p I .
\]
Inequality \eqref{eq:markov} has many applications in analysis, particularly
in approximation theory. For polynomial approximation and extension of
smooth functions see \cite{Levenberg,PP86,P22,PlesExt}. For inequalities on
polynomial and holomorphic images see \cite{P16,P21}.

\begin{definition}[The Siciak extremal function, also called the (pluricomplex) Green function]\label{def:siciak}
For $z\in\C^N$ put
\begin{equation}\label{eq:VK}
 V_K(z)=\sup\biggl\{\frac{\log|P(z)|}{n} :
              n\in\mathbb N,\ P\in\Pn_n,\ \norm{P}{K}\le1\biggr\}.
\end{equation}
Here $\log 0=-\infty$. The upper semicontinuous regularization of $V_K$ is
\[
 V_K^*(z)=\limsup_{w\to z}V_K(w).
\]
\end{definition}

The Siciak extremal function is one of the principal tools of
pluripotential theory \cite{DemTr,Klimek,Levenberg,Lev12,Siciak62,Sic81}.
The constant polynomial $1$ shows that $V_K\ge0$, and \eqref{eq:VK}
gives $V_K=0$ on $K$. The Lelong class $\mathcal L(\C^N)$ consists of the
plurisubharmonic functions $u$ on $\C^N$ for which there is a constant $C_u$
such that
\[
                   u(z)\le\log(1+|z|)+C_u \qquad(z\in\C^N).
\]
By Siciak's extremal theorem \cite[Theorem~4.12]{Sic81}, for every
$z\in\C^N$ we have
\[
 V_K(z)=\sup\bigl\{u(z) : u\in\mathcal L(\C^N),\
                    u(x)\le0\text{ for every }x\in K\bigr\}.
\]
Suppose that $N=1$ and that $K$ is nonpolar, that is, has positive
logarithmic capacity. Let $K_\infty$ be the component of
$\bigl(\C\cup\{\infty\}\bigr)\setminus K$ containing $\infty$.
Then $V_K$ is harmonic on $\C\setminus K$. At each finite point of
$K_\infty$ it agrees with the Green function of $K_\infty$ with pole at
infinity \cite[items~7.1 and~7.2]{Sic81}.

\begin{definition}[Pluripolarity]\label{def:polar}
A set $E\subset\C^N$ is pluripolar if there is a plurisubharmonic
function $u$ on $\C^N$, with $u\not\equiv-\infty$, such that
$u(x)=-\infty$ for every $x\in E$.
\end{definition}

A compact set $K$ is pluripolar if and only if $V_K^*\equiv+\infty$
\cite[Corollary~3.9 and Theorem~3.10]{Sic81}. If $K$ is nonpluripolar, then
$V_K^*$ is a locally bounded plurisubharmonic function in $\mathcal L(\C^N)$.
Its complex Monge--Ampère measure is therefore defined in the sense of
Bedford and Taylor \cite[p.~28]{Lev12}. This measure is supported on $K$
\cite[p.~33]{Lev12}, that is,
\[
               \bigl(dd^c V_K^*\bigr)^N=0\qquad\text{on }\C^N\setminus K.
\]

\begin{definition}[$L$-regularity]\label{def:regular}
The set $K$ is $L$-regular if $V_K$ is finite and continuous on $\C^N$.
\end{definition}

By Siciak \cite[Proposition~6.1]{Sic81}, $K$ is $L$-regular if and only if
$V_K$ is continuous at every $a\in K$, that is,
\[
               \lim_{\substack{z\to a\\ z\in\C^N}} V_K(z)=0
               \qquad(a\in K).
\]

\begin{definition}[The Hölder continuity property]\label{def:hcp}
For $0<\gamma\le1$, the set $K$ has the Hölder continuity property
with exponent $\gamma$, written $\mathrm{HCP}(\gamma)$, if there are
$B,r_0>0$ such that
\begin{equation}\label{eq:hcp}
 V_K(z)\le B\dist(z,K)^\gamma\qquad\bigl(z\in K_{r_0}\bigr).
\end{equation}
\end{definition}

Estimate \eqref{eq:hcp} is equivalent to Hölder continuity of $V_K$ with
exponent $\gamma$ on $\C^N$. Indeed, \eqref{eq:hcp} gives continuity at
each point of $K$. Hence $K$ is $L$-regular and $V_K\in\mathcal L(\C^N)$.
We use Błocki's translation argument in Siciak
\cite[Proposition~3.5]{Sic97several}. Normalizing $P(\cdot+h)$ in
\eqref{eq:VK} gives
\[
 V_K(z+h)-V_K(z)\le \max_{a\in K}V_K(a+h)
 \qquad\bigl(z,h\in\C^N\bigr).
\]
For $|h|\le r_0$, the right-hand side is at most $B|h|^\gamma$ by
\eqref{eq:hcp}. For $|h|>r_0$, it is at most $C+\log(1+|h|)$ by the
growth of $V_K$ and compactness of $K$, hence at most $C'|h|^\gamma$.
Interchanging $z$ and $z+h$ gives the Hölder estimate. The converse follows
by taking a point of $K$ nearest to $z$.

\subsection*{The main theorem}

\begin{theorem}\label{thm:main}
A nonempty compact set $K\subset\C^N$ satisfies Markov's inequality
for some $m$ if and only if it has the Hölder continuity property for
some $\gamma$. More precisely, \eqref{eq:hcp} implies \eqref{eq:markov}
with $m=1/\gamma$. Conversely, if \eqref{eq:markov} holds with constants
$A>0$ and $m\ge1$, then for every $0<r\le1/2$ and every $z\in K_r$,
\begin{equation}\label{eq:modulus}
 V_K(z)\le C r^{1/m}\biggl(\log\frac{e}{r}\biggr)^{1-1/m},
\end{equation}
where $C$ depends only on $N$, $A$, $m$ and $\operatorname{diam}K$.
If $N=1$, then for every $z\in\C$,
\begin{equation}\label{eq:plane}
 V_K(z)\le \frac{1+\log 8}{1-2^{-1/m}}\,
                    \bigl(A\dist(z,K)\bigr)^{1/m}.
\end{equation}
\end{theorem}

For every $0<\gamma<1/m$ the right-hand side of \eqref{eq:modulus} is
$O(r^\gamma)$ as $r\to0^+$.
Taking $r=\dist(z,K)>0$ gives \eqref{eq:hcp} near $K$. In particular, $K$ is
$L$-regular and nonpluripolar.

\begin{corollary}
Let $m_*(K)$ be the infimum of the admissible Markov exponents and let
$\gamma^*(K)$ be the supremum of the admissible Hölder exponents. For a
Markov set $K\subset \C^N$ we have
\[
                         \gamma^*(K)=\frac{1}{m_*(K)}.
\]
\end{corollary}

It was observed in \cite[p.~637]{BBM} that all examples then known
suggested this identity. Indeed, for every admissible
$\gamma$ and every admissible $m$, the two implications in
Theorem~\ref{thm:main} give
\[
 m_*(K)\le 1/\gamma,\qquad \gamma^*(K)\ge 1/m.
\]
If $N>1$ and \eqref{eq:markov} holds with $m=m_*(K)>1$, our proof
leaves open whether $K$ has $\mathrm{HCP}(1/m)$.

\begin{proof}[Proof that \eqref{eq:hcp} implies \eqref{eq:markov}]
We recall the standard argument for completeness.
Fix $x\in K$, $n\ge1$ and $P\in\Pn_n$ with $\norm P K\le1$.
Put $\rho=r_0n^{-1/\gamma}$. For each coordinate vector $e_j$,
\eqref{eq:VK} and \eqref{eq:hcp} give
\[
 |P(x+\zeta e_j)|\le \exp(nB\rho^\gamma)=\exp(Br_0^\gamma)
 \qquad(|\zeta|=\rho).
\]
Cauchy's estimate gives
\[
 |\partial_jP(x)|\le\frac{\exp(Br_0^\gamma)}{\rho},
 \qquad
 |\nabla P(x)|\le \frac{\sqrt N}{r_0}\exp(Br_0^\gamma)n^{1/\gamma}.
\]
Scaling proves \eqref{eq:markov}. If $\norm P K=0$, apply the same bound
to every scalar multiple of $P$.
\end{proof}

\section{A comparison of two neighborhoods}\label{sec:tubes}

\begin{lemma}\label{lem:tube}
For every $t>0$, the function $V_{K_t}$ is finite and plurisubharmonic on
$\C^N$. Moreover, it satisfies the following estimates.
\begin{enumerate}
\item[\textup{(i)}] For every $a\in K$ and $z\in\C^N$,
\[
 0\le V_{K_t}(z)\le\max\biggl\{0,\log\frac{|z-a|}{t}\biggr\}.
\]
In particular, taking $K=\{a\}$, we obtain for every $n\in\mathbb N$,
$P\in\Pn_n$ and $z\in\C^N$,
\begin{equation}\label{eq:ball}
 |P(z)|\le \max\biggl\{1,\frac{|z-a|^n}{t^n}\biggr\}
                  \,\norm P{\overline{B(a,t)}} .
\end{equation}
\item[\textup{(ii)}] For every $z,h\in\C^N$,
\begin{equation}\label{eq:lip}
 \bigl|V_{K_t}(z+h)-V_{K_t}(z)\bigr|
        \le \log\biggl(1+\frac{|h|}{t}\biggr).
\end{equation}
\end{enumerate}
\end{lemma}
\begin{proof}
Part~\textup{(i)} follows from Siciak's formula for the extremal function
of a ball \cite[p.~179, item~2.6]{Sic81} and
$\overline{B(a,t)}\subset K_t$.

For \textup{(ii)} we use Błocki's translation argument
(cf. \cite[Proposition~3.5]{Sic97several}). For $x\in K_t$ choose $a\in K$
with $|x-a|\le t$. Applying \eqref{eq:ball} at $x+h$ and taking the
supremum over $x\in K_t$, we obtain
\[
 \norm{P(\cdot+h)}{K_t}
 \le \biggl(1+\frac{|h|}{t}\biggr)^n\norm P{K_t}.
\]
Normalize the translated polynomial in \eqref{eq:VK} to get
\[
 V_{K_t}(z+h)\le V_{K_t}(z)+\log\biggl(1+\frac{|h|}{t}\biggr).
\]
Interchanging $z$ and $z+h$ gives \eqref{eq:lip}.
By \textup{(i)} and \eqref{eq:lip}, $V_{K_t}$ is finite and continuous.
Thus $V_{K_t}=V_{K_t}^*$. Since this regularization is finite, the
remark after Definition~\ref{def:polar} shows that $K_t$ is
nonpluripolar and $V_{K_t}$ is plurisubharmonic.
\end{proof}

For $t>0$ put
\[
                  b_K(t)=\sup_{z\in K_t}V_{K_{t/2}}(z).
\]
This quantity measures the growth of the extremal function between
two neighborhoods of $K$.

\begin{lemma}\label{lem:sum}
For every $r>0$ and $z\in\C^N$,
\begin{equation}\label{eq:sum}
 V_K(z)\le V_{K_r}(z)+\sum_{j=0}^{\infty}b_K(2^{-j}r).
\end{equation}
\end{lemma}
\begin{proof}
If $P\in\Pn_n$ and $\norm P{K_{t/2}}\le1$, then $\norm P{K_t}\le e^{nb_K(t)}$.
Thus for every $z\in\C^N$,
\[
 V_{K_{t/2}}(z)\le V_{K_t}(z)+b_K(t).
\]
For $q\ge1$, apply this inequality with $t=2^{-j}r$, $0\le j<q$, and add. We obtain
\[
 V_{K_{2^{\raisebox{1.4pt}{$\scriptscriptstyle -q$}}\,r}}(z)\le V_{K_r}(z)+\sum_{j=0}^{q-1}b_K(2^{-j}r).
\]
As $s\downarrow0$, $V_{K_s}$ increases to $V_K$. Indeed,
$\norm P{K_s}\downarrow\norm P K$, so every $P$ of norm at most one on
$K$ is the limit of the polynomials
\[
                 \frac{P}{\max\Bigl\{1,\norm P{K_s}\Bigr\}},
\]
each of norm at most one on $K_s$. Taking the supremum in \eqref{eq:VK}
and letting $q\to\infty$ proves \eqref{eq:sum}.
\end{proof}

For $n\in\mathbb N$ let
\[
 M_n(K)=\sup\bigl\{\norm{\nabla P}{K} :
                      P\in\Pn_n,\ \norm P K\le1\bigr\}.
\]
If $M_n(K)<\infty$, then
\[
 \norm{\nabla P}{K}\le M_n(K)\norm P K
                    \qquad\bigl(P\in\Pn_n\bigr).
\]

\begin{lemma}\label{lem:taylor}
If $M_n(K)<\infty$, then for every $t\ge0$ and $P\in\Pn_n$,
\begin{equation}\label{eq:taylor}
                  \norm P{K_t}\le e^{tM_n(K)}\norm P K .
\end{equation}
\end{lemma}
\begin{proof}
Fix $z\in K_t$, choose $a\in K$ nearest to $z$ and put $h=z-a$.
For $Q\in\Pn_n$, Cauchy--Schwarz and the definition of $M_n(K)$ give
\[
 \norm{D_hQ}K\le |h|\,\norm{\nabla Q}K
             \le |h|M_n(K)\norm QK .
\]
Since $D_hQ\in\Pn_n$, iteration gives
\[
 \norm{D_h^jP}K\le\bigl(|h|M_n(K)\bigr)^j\norm PK
 \qquad(1\le j\le n).
\]
Put $g(\zeta)=P(a+\zeta h)$ for $\zeta\in\C$. By the chain rule,
$g^{(j)}(0)=D_h^jP(a)$ for $1\le j\le n$. Taylor's formula gives
\[
 P(z)=g(1)=P(a)+\sum_{j=1}^n\frac{D_h^jP(a)}{j!}.
\]
Taking absolute values and using the preceding bound and $|h|\le t$, we obtain
\[
 |P(z)|\le \Biggl(1+\sum_{j=1}^n\frac{\bigl(|h|M_n(K)\bigr)^j}{j!}\Biggr)
                \norm PK
        \le e^{tM_n(K)}\norm PK .\qedhere
\]
\end{proof}

By \eqref{eq:taylor}, every Markov set determines polynomials. Indeed, if $P$ vanishes on $K$, it vanishes on every $K_t$.

\begin{corollary}\label{cor:linear}
Suppose that $A>0$ satisfies
\[
 \norm{\nabla P}K\le An\norm P K
 \qquad\bigl(n\in\mathbb N,\ P\in\Pn_n\bigr).
\]
Then for every $z,w\in\C^N$,
\[
 V_K(z)\le\log\bigl(1+A\dist(z,K)\bigr),\qquad
 |V_K(z)-V_K(w)|\le\log\bigl(1+A|z-w|\bigr).
\]
\end{corollary}
\begin{proof}
Fix a unit vector $v\in\C^N$. We show by induction on $n$ that
\[
 |P(a+rv)|\le(1+Ar)^n\norm P K
 \qquad\bigl(P\in\Pn_n,\ a\in K,\ r\ge0\bigr).
\]
For $n=0$ the assertion is clear. Suppose it holds for $n-1$, where
$n\ge1$. Put $Q=D_vP$. Then $Q\in\Pn_{n-1}$ and
\[
 \norm Q K\le\norm{\nabla P}K\le An\norm P K .
\]
For every $t\in\mathbb R$,
\[
 \frac{d}{dt}P(a+tv)=Q(a+tv).
\]
Integrating from $0$ to $r$ and applying the induction hypothesis to $Q$ gives
\begin{align*}
 |P(a+rv)| &\le \norm P K+\int_0^r |Q(a+tv)|\,dt\\
 &\le \norm P K\biggl(1+An\int_0^r(1+At)^{n-1}\,dt\biggr)
       =(1+Ar)^n\norm P K .
\end{align*}
Thus for every $h\in\C^N$,
\[
 \norm{P(\cdot+h)}K\le(1+A|h|)^n\norm PK .
\]
Normalizing this translate in \eqref{eq:VK}, as in the proof of
Lemma~\ref{lem:tube}\,\textup{(ii)}, gives
\[
 V_K(w+h)\le V_K(w)+\log\bigl(1+A|h|\bigr).
\]
Taking $w\in K$ first shows that $V_K$ is finite on $\C^N$.
For arbitrary $z,w\in\C^N$, put $h=z-w$ and interchange $z$ and $w$.
This gives the second estimate. Taking $w\in K$ nearest to $z$ gives
the first estimate.
\end{proof}

Corollary~\ref{cor:linear} is the well-known case $m=1$.
For a planar set whose unbounded complementary component is regular, see
Toókos and Totik \cite{TT05}. For proofs of Lipschitz continuity using Taylor's formula, see
\cite[Theorem~5.3]{BC11} and \cite[pp.~637--638]{BBM}.
The bound is attained on a closed disk of radius
$R>0$ in $\C$, with $A=1/R$.

For $m>1$, the same induction works up to degree $n$ with $An^{m-1}$
in place of $A$.
After taking the $n$th root, the factor $1+Ar n^{m-1}$ is
unbounded in $n$ at every fixed distance $r=\dist(z,K)>0$.
Our aim in Theorem~\ref{thm:main} is to prove that \eqref{eq:markov}
still implies Hölder continuity of $V_K$ with every exponent
$0<\gamma<1/m$.

\section{The estimate in several variables}\label{sec:general}

To prove the next lemma, we apply Hörmander's extension theorem
\cite[Theorem~4.4.3]{Hormander}.

\begin{lemma}\label{lem:point}
Let $s>0$, $k\in\mathbb N$, $y\in\C^N$ and $W=V_{K_s}$. There is
$Q\in\Pn_{k+3N}$ such that $Q(y)=e^{kW(y)}$ and, for every $w\in\C^N$,
\[
 |Q(w)|\le C_N\biggl(1+\frac{2k|w-y|}{s}\biggr)^{3N}e^{kW(w)}.
\]
The constant $C_N\ge1$ depends only on $N$.
\end{lemma}
\begin{proof}
Put $\delta=s/(2k)$ and $\varphi(\xi)=2kW(y+\delta\xi)$. By
Lemma~\ref{lem:tube}, $\varphi$ is plurisubharmonic and, for $\xi,\eta\in\C^N$,
\[
 \bigl|\varphi(\xi)-\varphi(\eta)\bigr|
   \le 2k\log\biggl(1+\frac{|\xi-\eta|}{2k}\biggr)\le |\xi-\eta|.
\]
Thus Hörmander's condition (4.4.7) holds with oscillation
constant equal to one \cite[pp.~95--96]{Hormander}.
In his Theorem~4.4.3, formula~(4.4.9), take the subspace $\Sigma=\{0\}$
of complex codimension $N$ and prescribe $u(0)=e^{\varphi(0)/2}$.
Integration over this subspace gives
\[
           \int_\Sigma |u|^2e^{-\varphi}\,d\sigma
                     =|u(0)|^2e^{-\varphi(0)}=1.
\]
Consequently, we obtain an entire function $G$ with $G(0)=u(0)$ such that
\[
 \int_{\C^N}|G(\eta)|^2e^{-\varphi(\eta)}
           \bigl(1+|\eta|^2\bigr)^{-3N}\,dV(\eta)\le(9\pi e)^N ,
\]
where $dV$ denotes Lebesgue measure.

Fix $\xi\in\C^N$. For $\eta\in B(\xi,1)$ we have
\[
 \varphi(\eta)\le\varphi(\xi)+1,\qquad
 1+|\eta|^2\le 1+(|\xi|+1)^2\le2(1+|\xi|)^2 .
\]
The mean value inequality for $|G|^2$ on the ball $B(\xi,1)$ gives
\[
 |G(\xi)|^2\le \frac{1}{\operatorname{vol}\bigl(B(0,1)\bigr)}
                 \int_{B(\xi,1)}|G(\eta)|^2\,dV(\eta)
          \le C_N^2(1+|\xi|)^{6N}e^{\varphi(\xi)} .
\]
Define
\[
 Q(w)=G\biggl(\frac{w-y}{\delta}\biggr)\qquad(w\in\C^N).
\]
Then $Q(y)=e^{kW(y)}$, and the preceding inequality gives the asserted
bound for $Q$. Fix $a\in K$.
By Lemma~\ref{lem:tube}\,\textup{(i)},
\[
 e^{W(w)}\le\max\biggl\{1,\frac{|w-a|}{s}\biggr\},
 \qquad |Q(w)|\le M(1+|w|)^{k+3N},
\]
where $M$ is independent of $w$. The polynomial-growth form of Liouville's
theorem gives $Q\in\Pn_{k+3N}$.
\end{proof}

\begin{proof}[Proof of \eqref{eq:modulus}]
For $m=1$ the assertion follows from Corollary~\ref{cor:linear} and
$\log(1+Ar)\le Ar$. Suppose $m>1$. Fix $0<t\le1/2$, put $s=t/2$ and
take $y\in K_t$. Write $D=\operatorname{diam}K$. For $x\in K$ we have
$V_{K_s}(x)=0$ and $|x-y|\le D+t$. The polynomial of
Lemma~\ref{lem:point} therefore satisfies
\[
          \norm Q K\le C_N\biggl(1+\frac{4k(D+t)}{t}\biggr)^{3N}.
\]
Since $Q\in\Pn_{k+3N}$, formulas \eqref{eq:taylor} and \eqref{eq:markov} give
\[
 e^{kV_{K_s}(y)}=|Q(y)|\le e^{At(k+3N)^m}\norm Q K .
\]
Taking logarithms and the supremum over $y\in K_t$, we obtain
\[
 b_K(t)\le \frac{At(k+3N)^m}{k}
       +\frac{1}{k}\biggl[\log C_N+
                     3N\log\biggl(1+\frac{4k(D+t)}{t}\biggr)\biggr].
\]
Choose
\[
 \ell(t)=\log\frac{e}{t},\qquad
 k=\left\lceil\biggl(\frac{\ell(t)}{t}\biggr)^{1/m}\right\rceil .
\]
The inequalities
\[
 k+3N\le(1+3N)k,\qquad k\le2\biggl(\frac{\ell(t)}{t}\biggr)^{1/m}
\]
give the bound $C_{N,A,m}t^{1/m}\ell(t)^{1-1/m}$ for the first term.
For the second term, use $t\le1/2$ and $k\ge1$ to obtain
\[
 \log\biggl(1+\frac{4k(D+t)}{t}\biggr)
    \le \log\bigl(1+4(D+1)\bigr)+\log k+\log(1/t)
    \le C_{m,D}\ell(t).
\]
The lower bound $k\ge\bigl(\ell(t)/t\bigr)^{1/m}$ now gives
$C_{N,m,D}t^{1/m}\ell(t)^{1-1/m}$ as a bound for the second term.

Now fix $0<r\le1/2$ and $z\in K_r$. Then $V_{K_r}(z)=0$. Moreover,
\[
 \ell(2^{-j}r)=\ell(r)+j\log 2\le(1+j)\ell(r).
\]
Substitution in Lemma~\ref{lem:sum} gives
\[
 V_K(z)\le C_{N,A,m,D}\,r^{1/m}\ell(r)^{1-1/m}
             \sum_{j=0}^{\infty}2^{-j/m}(1+j)^{1-1/m}.
\]
The series converges, which proves \eqref{eq:modulus}.
\end{proof}

\section{An elementary proof on the complex line}\label{sec:plane}

\begin{proof}[Proof of \eqref{eq:plane}]
Fix $t>0$ and $n\in\mathbb N_0$. Put $s=t/2$. We first prove that
\begin{equation}\label{eq:planar-step}
 (n+1)V_{K_s}(\xi)\le Atn^m+\log 8
                    \qquad\bigl(\xi\in K_t\bigr).
\end{equation}
For $\xi\in K_s$ this follows from $V_{K_s}(\xi)=0$. Fix
$\xi\in K_t\setminus K_s$ and choose $a\in K$ nearest to $\xi$.
We claim that there are finitely many $x_i\in K$ and $c_i\in\C$,
depending on $t,n,\xi$, such that
\begin{equation}\label{eq:moments}
 F(\xi)=\sum_i c_iF(x_i)\quad\bigl(F\in\Pn_n\bigr),\qquad
 \sum_i|c_i| = e^{Atn^m}.
\end{equation}
Put $X(\zeta)=(1,\zeta,\ldots,\zeta^n)$ for $\zeta\in\C$. Let $H$ be the
convex hull in $\C^{n+1}$ of the compact set
\[
 \Bigl\{e^{i\theta}X(x) : x\in K,\ 0\le\theta\le2\pi\Bigr\}.
\]
The set $H$ is compact.
We show that $y=e^{-Atn^m}X(\xi)$ belongs to $H$. Suppose otherwise and choose
$q\in H$ nearest to $y$. For $v=(v_0,\ldots,v_n)\in\C^{n+1}$ put
\[
 L(v)=\operatorname{Re}\sum_{j=0}^n \overline{y_j-q_j}\,v_j .
\]
For every $v\in H$, the segment $q+\tau(v-q)$, $0\le\tau\le1$, lies in $H$.
Thus $\tau\mapsto |y-q-\tau(v-q)|^2$ has a minimum at $\tau=0$ on $[0,1]$.
Its right derivative there is $-2\bigl(L(v)-L(q)\bigr)$, so $L(v)\le L(q)$.
Also $L(y)-L(q)=|y-q|^2>0$.

Consider the polynomial $\phi(\zeta)=\sum_{j=0}^n\overline{y_j-q_j}\,\zeta^j$.
For every $x\in K$ we have
\[
 |\phi(x)|=\max_{0\le\theta\le2\pi}\operatorname{Re}\bigl(e^{i\theta}\phi(x)\bigr)
       =\max_{0\le\theta\le2\pi}L\bigl(e^{i\theta}X(x)\bigr)\le L(q).
\]
Since $y_j=e^{-Atn^m}\xi^j$, it follows that
\[
 \norm\phi K\le L(q)<L(y)=e^{-Atn^m}\operatorname{Re}\phi(\xi)
       \le e^{-Atn^m}|\phi(\xi)|\le\norm\phi K .
\]
For $n\ge1$ the last inequality follows from \eqref{eq:taylor} and
\eqref{eq:markov}. For $n=0$ it follows from the fact that $\phi$ is constant.
This contradiction proves that $y\in H$. Hence
\[
 e^{-Atn^m}X(\xi)=\sum_i\lambda_i e^{i\theta_i}X(x_i),
 \qquad \lambda_i\ge0,\qquad \sum_i\lambda_i=1 .
\]
Put $c_i=e^{Atn^m}\lambda_i e^{i\theta_i}$. Then $\sum_i|c_i|=e^{Atn^m}$.
Comparing coordinates gives the first identity in \eqref{eq:moments}
for $1,x,\ldots,x^n$, and hence for every $F\in\Pn_n$ by linearity.

For $w\in\C\setminus K$ set
\begin{equation}\label{eq:planar-g}
 g(w)=\sum_i c_i\frac{\xi-x_i}{w-x_i},\qquad h(w)=\frac{\xi-a}{w-a}\,g(w).
\end{equation}
Then
\[
 h(\xi)=g(\xi)=\sum_i c_i=1.
\]
For $0\le j<n$, apply \eqref{eq:moments} to $(\xi-x)x^j$ as a polynomial
in $x$. We get
\[
 \sum_i c_i(\xi-x_i)x_i^j=0.
\]
Since $(w^n-x^n)/(w-x)$ is a polynomial
in $x$ of degree less than $n$ for $n\ge1$, these identities give
\[
 w^n g(w)=\sum_i \frac{c_i(\xi-x_i)x_i^n}{w-x_i}.
\]
For $n=0$ this is the definition of $g$. Thus $g(w)=O(|w|^{-n-1})$
and $h(w)=O(|w|^{-n-2})$ at infinity.

For $w\in\partial K_s$ and $x\in K$ we have $|w-a|,|w-x|\ge s$.
Since $|\xi-a|\le t=2s$ and
\[
 |\xi-x|\le |w-x|+|w-a|+|\xi-a|,
\]
we obtain
\[
 \frac{|\xi-a|}{|w-a|}\frac{|\xi-x|}{|w-x|}
 \le 2\frac{|\xi-a|}{s}+\frac{|\xi-a|^2}{s^2}\le8.
\]
Thus $|h(w)|\le8e^{Atn^m}$ on $\partial K_s$. Fix $d\in\mathbb N$ and
$P\in\Pn_d$ with $\norm P{K_s}\le1$. The function $P^{n+1}h^d$ is
holomorphic on $\C\setminus K_s$ and tends to zero at infinity.
On $\partial K_s$ its modulus is at most $\bigl(8e^{Atn^m}\bigr)^d$.
Choose a disk whose interior contains $K_s$ and $\xi$, and on whose
boundary the same bound holds. The maximum principle gives this bound
on the component of the disk minus $K_s$ containing $\xi$. Since $h(\xi)=1$,
\[
 |P(\xi)|^{n+1}\le\bigl(8e^{Atn^m}\bigr)^d .
\]
Taking logarithms, dividing by $d$ and taking the supremum over all such
$d$ and $P$ gives \eqref{eq:planar-step}.

To finish the proof, put $n=\bigl\lfloor(At)^{-1/m}\bigr\rfloor$. Then
\[
 Atn^m\le1,\qquad \frac{1}{n+1}\le(At)^{1/m}.
\]
By \eqref{eq:planar-step}, for every $\xi\in K_t$ we have
\[
 V_{K_s}(\xi)\le(1+\log 8)(At)^{1/m}.
\]
Taking the supremum over $\xi\in K_t$ gives
\begin{equation}\label{eq:planar-b}
 b_K(t)\le(1+\log 8)(At)^{1/m}.
\end{equation}
Fix $z\in\C$. If $z\in K$, both sides of \eqref{eq:plane} vanish.
Otherwise put $r=\dist(z,K)$. Then $V_{K_r}(z)=0$, and
\eqref{eq:sum} and \eqref{eq:planar-b} give
\[
 V_K(z)\le (1+\log 8)(Ar)^{1/m}\sum_{j=0}^{\infty}2^{-j/m}
          =\frac{1+\log 8}{1-2^{-1/m}}(Ar)^{1/m}.\qedhere
\]
\end{proof}

\begin{corollary}[Products]\label{cor:product}
Suppose $K=E_1\times\cdots\times E_N$, with each $E_j\subset\C$
nonempty and compact. If $K$ satisfies \eqref{eq:markov} with constants
$A>0$ and $m\ge1$, then for every $z\in\C^N$,
\[
 V_K(z)\le C_m\bigl(A\dist(z,K)\bigr)^{1/m},\qquad
 C_m=\frac{1+\log8}{1-2^{-1/m}}.
\]
\end{corollary}
\begin{proof}
Testing polynomials in one coordinate gives \eqref{eq:markov} on
each $E_j$ with the same $A,m$. By Siciak's product formula
\cite[Proposition~5.9]{Sic81} and \eqref{eq:plane},
\[
 V_K(z)=\max_j V_{E_j}(z_j)
       \le C_m\max_j\bigl(A\dist(z_j,E_j)\bigr)^{1/m}
       \le C_m\bigl(A\dist(z,K)\bigr)^{1/m}.\qedhere
\]
\end{proof}

\section*{Acknowledgements}
During the preparation of this manuscript, the authors used ChatGPT (OpenAI) for language editing and to explore possible simplifications of certain mathematical arguments. The authors critically examined the suggestions provided by the tool, substantially modifying or entirely replacing some of them. All mathematical statements and proofs were independently verified by the authors, who take full responsibility for the content of this work.

\end{document}